\documentclass[twoside,10pt,letter]{article}
\usepackage[margin=1.2in]{geometry}

\usepackage{times}
\usepackage{amsmath}
\usepackage{amsfonts}
\usepackage{amssymb}
\usepackage{amsthm}
\usepackage{makeidx} 
\usepackage{graphicx}
\usepackage{mathrsfs}
\usepackage{tabularx}
\usepackage{caption}
\usepackage{bbold}
\usepackage{color}
\usepackage{fancybox}
\usepackage{verbatim}
\usepackage{hyperref}
\usepackage{tikz}
\usepackage{float}
\usepackage{url}
\usetikzlibrary{arrows,automata}
\usetikzlibrary{trees}
\usetikzlibrary{shadows}
\usepackage[framemethod=TikZ]{mdframed}
\usepackage{mdframed}
\mdfdefinestyle{theorem}{%
	linecolor=black,
	outerlinewidth=.5pt,
	roundcorner=5pt,
	innertopmargin=.5\baselineskip,
	innerbottommargin=.75\baselineskip,
	innerrightmargin=20pt,
	innerleftmargin=20pt,
	backgroundcolor=green!20,
	shadow=true,
	shadowcolor=green!20}
\mdfdefinestyle{algo}{%
	linecolor=black,
	outerlinewidth=.5pt,
	roundcorner=5pt,
	innertopmargin=.5\baselineskip,
	innerbottommargin=.75\baselineskip,
	innerrightmargin=10pt,
	innerleftmargin=10pt,
	backgroundcolor=yellow!20,
	shadow=true,
	shadowcolor=yellow!20}
\mdfdefinestyle{example}{%
	linecolor=black,
	outerlinewidth=.5pt,
	leftline = false,
	rightline = false,
	innertopmargin=.5\baselineskip,
	innerbottommargin=.75\baselineskip,
	innerrightmargin=10pt,
	innerleftmargin=10pt,
	backgroundcolor=white}

\def\bP{\mathbb{P}}
\def\bR{\mathbb{R}}

\def\b1{\mathbb{1}}

\def\e1{\mathsf{1}}

\def\tT{\mathtt{T}}

\def\cU{{\mathcal{U}}}

\def\e1{\mathsf{1}}

\def\of0{(0)}

\def\tr{\mathrm{tr}}

\def\bf0{\mathbf{0}}
\def\cp1{\mathbb{CP}^1}

\def\arg{\mathrm{arg}}

\theoremstyle{definition}

\theoremstyle{definition}

\begin{document}

\title{\textbf{Differential Evolution for Grassmann Manifold\\ Optimization: A Projection Approach}}
\author{\textbf{Andrew Lesniewski}\\
Department of Mathematics\\
Baruch College\\
One Bernard Baruch Way\\
New York, NY 10010\\
USA}
\maketitle

\begin{abstract}
We propose a novel evolutionary algorithm for optimizing real-valued objective functions defined on the Grassmann manifold $\mathrm{Gr}(k,n)$, the space of all $k$-dimensional linear subspaces of $\bR^n$. While existing optimization techniques on $\mathrm{Gr}(k,n)$ predominantly rely on first- or second-order Riemannian methods, these inherently local methods often struggle with nonconvex or multimodal landscapes. To address this limitation, we adapt the Differential Evolution algorithm - a global, population based optimization method - to operate effectively on the Grassmannian. Our approach incorporates adaptive control parameter schemes, and introduces a projection mechanism that maps trial vectors onto the manifold via QR decomposition. The resulting algorithm maintains feasibility with respect to the manifold structure while enabling exploration beyond local neighborhoods. This framework provides a flexible and geometry-aware alternative to classical Riemannian optimization methods and is well-suited to applications in machine learning, signal processing, and low-rank matrix recovery where subspace representations play a central role. We test the methodology on a number of examples of optimization problems on Grassmann manifolds.
\end{abstract}

\section{Introduction}

The (real) Grassmann manifold, or Grassmannian, $\mathrm{Gr}(k,n)$ is the set of all $k$-dimensional linear subspaces of $\bR^n$. Formally,
\begin{equation}
\mathrm{Gr}(k,n)=\{ V \subseteq \bR^n \mid \dim V = k\}.
\end{equation}
A special case of the Grassmannian is a (real) projective space $\bP\bR^{n-1}\cong\mathrm{Gr}(1,n)$, which corresponds to the space of lines through the origin in $\bR^n$.

The Grassmannian $\mathrm{Gr}(k,n)$ possesses a rich geometric structure: it is a smooth, compact manifold of dimension
\begin{equation}
\dim \mathrm{Gr}(k,n) = k(n - k),
\end{equation}
and can be naturally endowed with a Riemannian structure. This enables the use of differential geometric tools for optimization and learning on the manifold \cite{AMS08}, \cite{B23}.

A common representation of the Grassmannian is via its construction as a quotient of the \textit{Stiefel manifold}. The Stiefel manifold $\mathrm{V}_k(\bR^n)$ is the space of all orthonormal $k$-frames in $\bR^n$, i.e., matrices $Q\in \mathrm{Mat}_{n,k}(\bR)$ satisfying:
\begin{equation}\label{orthcond}
Q^\tT Q = I_k.
\end{equation}
The Grassmannian is then obtained by identifying frames that differ by an orthogonal transformation in $O(k)$, yielding the quotient:
\begin{equation}
\mathrm{Gr}(k,n) = \mathrm{V}_k(\bR^n) \big/ O(k),
\end{equation}
where the group $O(k)$ acts on the right by rotation,
\begin{equation}\label{groupact}
Q \sim Q R,
\end{equation}
for $R \in O(k)$.

Alternatively, any $k$-dimensional subspace of $\bR^n$ can be represented by a full-rank matrix $Y \in\mathrm{Mat}_{n,k}(\bR)$, whose columns span the subspace. Two matrices $Y_1, Y_2 \in \mathrm{Mat}_{n,k}(\bR)$ represent the same point in the Grassmannian if
\begin{equation*}
\operatorname{span}(Y_1) = \operatorname{span}(Y_2).
\end{equation*}
In particular, matrices $Q \in\mathrm{Mat}_{n,k}(\bR)$, obtained via QR decomposition of $Y$ and satisfying \eqref{orthcond}, serve as canonical representatives of points on $\mathrm{Gr}(k,n)$.

Grassmannians play a fundamental - albeit specialized - role in various modern applications across mathematics and applied sciences. In artificial intelligence and machine learning, they are particularly useful in settings where data is better modeled by subspaces rather than points, giving rise to geometry aware learning. A partial list of such applications includes: machine learning and computer vision (subspace-based classification, face recognition, action recognition), signal processing (subspace tracking, source separation), optimization (low-rank matrix recovery, dictionary learning on manifolds), and statistics (canonical correlation analysis, dimension reduction, multivariate regression). These applications benefit from the intrinsic geometry of the Grassmannian, enabling more robust, structured, and interpretable models when working with high-dimensional data that inherently varies along subspaces. For further details and representative developments, see \cite{FWLLC23}, \cite{EV13}, \cite{TVSC11}, \cite{ZZHH18}.

Optimization on Grassmann manifolds has been an area of active and growing research, owing to its relevance in applications where data or model parameters naturally reside on spaces of linear subspaces. Central to this body of work are first- and second-order optimization methods that leverage the underlying Riemannian geometry of the Grassmannian $\mathrm{Gr}(k,n)$ \cite{AMS08, B23}. These methods extend classical optimization techniques - such as line search, gradient descent, and Newton's method - to curved spaces by replacing Euclidean concepts with their manifold analogues (e.g., Riemannian gradients, exponential maps, parallel transport).

While these Riemannian approaches are elegant and often highly efficient in practice, they are inherently local in nature. That is, they rely on smoothness and local curvature information to make progress, and as such, they do not guarantee convergence to global optima of nonconvex (in the ambient space) objective functions. This limitation is particularly pronounced in high-dimensional, multi-modal landscapes where local methods may become trapped in suboptimal basins of attraction.

In this note, we propose a global, population based optimization algorithm tailored for functions defined on Grassmann manifolds. Our method is grounded in the Differential Evolution (DE) framework \cite{SP97}, a robust and widely used evolutionary algorithm for continuous optimization. To enhance adaptability and improve convergence behavior across varying objective landscapes, we incorporate the self-adaptive parameter control strategies introduced in \cite{BGBMZ06}. These mechanisms dynamically adjust key algorithmic parameters, such as mutation and crossover rates, thereby reducing the need for manual tuning and increasing robustness across diverse problem instances.

Crucially, we adapt the DE framework to respect the geometric constraints of the Grassmannian by introducing a projection step that ensures all candidate solutions remain on the manifold. This projection is implemented via QR decomposition and is seamlessly integrated into the mutation and crossover phases. The result is a flexible and geometry aware optimization algorithm capable of exploring the Grassmannian effectively while preserving the essential structure of the search space.

Our approach is particularly well suited for problems in which global search capabilities are essential, and where the solution is naturally expressed as a subspace such as in dimensionality reduction, subspace learning, and low rank modeling. The algorithm demonstrates strong empirical performance for problems of moderate size, especially when the ambient dimension $d=nk$ does not significantly exceed 100.

\section{The algorithm}

We consider the problem of minimizing a function $f(x)$, where $x \in\bR^d$, with $d$ representing the dimensionality of the search space. The DE algorithm is a powerful population-based optimization technique designed to explore this space efficiently. The fundamental idea behind DE is a structured approach to generating trial parameter vectors that iteratively refine potential solutions. 

The classical version of DE, introduced in \cite{SP97}, consists of three operations: mutation, crossover, and selection. The algorithm begins by initializing a population of $NP\geq 4$ candidate solutions $x_{i,G}$, $i=1,\ldots,NP$, each represented as a $d$-dimensional vector. The index $G$ denotes the current generation, with the initial population corresponding to $G=0$. 

In the mutation phase, new candidate solutions (mutant vectors) are generated by combining existing ones in a stochastic manner. The crossover phase then blends the mutant vectors with members of the current population to create trial vectors, introducing diversity into the population. Finally, the selection phase determines which vectors will progress to the next generation $G\to G+1$, favoring those that yield better objective function values. This evolutionary cycles continues until a predefined stopping criterion - such as a convergence threshold and a maximum number of generations - is satisfied.

The effectiveness of DE is highly dependent on two key control parameters: $F\in[0,2]$, the differential weight, which controls the magnitude of mutation, and $CR\in(0,1)$, the crossover probability, which influences the extent to which trial vectors inherit components from mutant vectors. These parameters directly influence the balance between exploration and exploitation, and thus play a critical role in the convergence behavior of the algorithm.

The three main operations of DE are defined as follows:

\begin{mdframed}[style=algo]
\textbf{Mutation}. For each target vector $x_{i,G}$, a mutant vector $v_{i,G+1}$ is generated according to
\begin{equation}
v_{i,G+1}=x_{r_1,G}+F(x_{r_2,G}-x_{r_3,G}), 
\end{equation}
where the indices $1\leq r_1,r_2, r_3\leq NP$ are randomly chosen, in such a way that the four indices $i, r_1,r_2, r_3$ are distinct.
\end{mdframed}
This ensures diversity in the mutation process by preventing self-selection. The scaling factor $F$ controls the amplitude of the differential variation and thereby the algorithm’s exploratory behavior. 

\begin{mdframed}[style=algo]
\textbf{Crossover}. To enhance diversity and facilitate exploration, the crossover step generates a trial vector $u_{i,G+1}$ by combining the mutant vector $v_{i,G+1}$ with the target vector $x_{i,G}$, $i=1,\ldots,NP$. Namely, we generate $d$ random numbers $\alpha^j$ from the uniform distribution $\cU(0,1)$, and $NP$ random integers $1\leq R_i\leq NP$. The target vector is mixed with the mutated vector, using the following scheme, to yield the trial vector
\begin{equation}
u_{i,G+1}=(u^1_{i,G+1},\ldots,u^d_{i,G+1}),
\end{equation}
where
\begin{equation}
u^j_{i,G+1}=
\begin{cases}
v^j_{i,G+1},&\text{ if }\alpha^j\leq CR\text{ or } j=R_i,\\
x_{i,G}&\text{ if }\alpha^j> CR\text{ and } j\neq R_i,
\end{cases}
\end{equation}
for $j=1,\ldots,d$.
\end{mdframed}
This mechanism guarantees that the trial vector differs from its parent in at least one dimension, thereby preserving diversity in the population and avoiding stagnation.

\begin{mdframed}[style=algo]
\textbf{Selection}. The selection phase employs a greedy approach to determine the next generation population:
\begin{equation}
x_{i,G+1}=
\begin{cases}
u_{i,G+1},&\text{ if } f(u_{i,G+1}) < f(x_{i,G}),\\
x_{i,G},&\text{ otherwise, }\\
\end{cases}
\end{equation}
for $i=1,\ldots,NP$. 
\end{mdframed}
This ensures that the population either improves or remains unchanged in each iteration, thus maintaining the algorithm's convergence properties.

While DE is simple and effective, its performance is highly sensitive to the choice of the control parameters $F$ and $CR$. Improper parameter settings can lead to premature convergence, stagnation, or excessive random search. To address this challenge, several adaptive and self-adaptive variants of DE have been developed, where parameters are dynamically tuned during the evolution process. In this work, we adopt the adaptation strategy introduced by \cite{BGBMZ06}, which has demonstrated strong empirical performance, especially in problems with moderate dimensionality (up to $d = 100$).

In this adaptive scheme, each individual target vector $x_{i,G}$ in the population is associated with its own evolving parameters $F_{i,G}$ and $CR_{i,G} $. These parameters are updated at each generation according to probabilistic rules, allowing the algorithm to self-tune and balance exploration and exploitation over time.

\begin{mdframed}[style=algo]
\textbf{Parameter Adaptation}. For each individual $i = 1,\ldots, NP$, the control parameters for the next generation are updated as follows:
\begin{equation}
\begin{split}
F_{i,G+1} &=
\begin{cases}
F_l + r_1 F_u, & \text{if } r_2 < \tau_1, \\
F_{i,G}, & \text{otherwise},
\end{cases} \\
CR_{i,G+1} &=
\begin{cases}
r_3, & \text{if } r_4 < \tau_2, \\
CR_{i,G}, & \text{otherwise},
\end{cases}
\end{split}
\end{equation}
where: $r_1, r_2, r_3, r_4 \sim \cU(0, 1)$ are independent random variables drawn from the uniform distribution, $F_l$ is the minimum value for $F$, $F_u$ is the upper range for perturbation (so $F_{i,G+1} \in [F_l, F_l + F_u]$), and $\tau_1$ and $\tau_2$ are the probabilities controlling how often the parameters are updated. A commonly used configuration is: $F_l=0.1$, $F_u=0.9$, $\tau_1=\tau_2=0.1$.
\end{mdframed}
This stochastic adaptation mechanism promotes behavioral diversity across the population, helping the algorithm escape local optima and enhancing its ability to explore complex, high-dimensional landscapes.

The final operation necessary to ensure the effectiveness of the algorithm when applied to the Grassmannian manifold $\mathrm{Gr}(k,n)$ - the space of all $k$-dimensional subspaces of $\bR^n$ is the projection operation. Since each individual in the population must represent a point on $\mathrm{Gr}(k,n)$, we structure the search space accordingly by setting the ambient dimension to $d=nk$, treating each vector as a flattened matrix in 
$\mathrm{Mat}_{n,k}(\bR)$.
\begin{mdframed}[style=algo]
\textbf{Projection}.
Given a population of vectors $y_{i, G}$, $i=1,\ldots,NP$, where $d=nk$, perform the following steps:
\begin{itemize}
\item[(i)] Reshape each vector into an $n\times k$ matrix $Y_{i,G}$
\item[(ii)] Perform the QR decomposition on $Y_{i,G}$, $Y_{i,G}=Q_{i,G}R_{i,G}$, where $Q_{i,G}\in\mathrm{Mat}_{n,k}(\bR)$ is an orthonormal matrix.
\item[(iii)] Take the transpose $Q_{i,G}^\tT$, flatten it into a vector in $\bR^d$, and use this as the projected vector representing a point on the Grassmannian.
\end{itemize}
\end{mdframed}
This ensures that each individual lies on $\mathrm{Gr}(k,n)$, preserving feasibility with respect to the underlying manifold structure.

The complete algorithm is formulated as follows.
\begin{mdframed}[style=algo]
\begin{equation}
\begin{aligned}
&pop,\ F,\ CR \gets \mathrm{Initialize}() \\
&\mathrm{Project}(pop) \\
&\textbf{while}\ \text{exit criteria not met} \\
&\qquad \mathrm{Mutate}(pop) \\
&\qquad \mathrm{Project}(pop) \\
&\qquad \mathrm{Crossover}(pop) \\
&\qquad \mathrm{Project}(pop) \\
&\qquad \mathrm{Select}(pop) \\
&\qquad F,\ CR \gets \mathrm{Adapt}() \\
&\textbf{return}\ pop
\end{aligned}
\end{equation}
\end{mdframed}

\section{Experiments}

We consider five nontrivial test functions defined on the Grassmannian $\mathrm{Gr}(k,n)$ to evaluate the performance of the proposed optimization algorithm. These functions are representative of objective landscapes that naturally arise in applications such as machine learning, signal processing, and statistical modeling, where subspace-based representations play a central role.

Each of our test functions is based on a fundamental functional form:
\begin{equation}\label{basicfct}
\phi(Q)=\tr(Q^\tT A Q),
\end{equation}
where $Q\in\mathrm{Mat}_{n,k}(\bR)$ is an orthonormal matrix and $A\in\mathrm{Mat}_n(\bR)$ is a fixed symmetric matrix. Due to the invariance of \eqref{basicfct} under the right action of the orthogonal group \eqref{groupact}, functions of this form depend only on the coset class $V$ (the $k$-dimensional subspace spanned by $Q$) of $Q$ in the Grassmannian. Accordingly, we will henceforth write $\phi(V)$ to denote its value on the Grassmannian.

By $\|\cdot\|_F$, we will denote the Frobenius norm. Note that, in particular,
\begin{equation}
\|Q^\tT P\|_F^2=\tr(Q^\tT PP^\tT Q),
\end{equation}
which is of the form \eqref{basicfct}. For $Q,P\in\mathrm{Gr}(k,n)$, $Q^\tT P=\{q_i^\tT p_j\}_{1\leq i,j\leq k}$, where $q_i$ and $p_j$ denote the columns of $Q$ and $P$, respectively. The squared Frobenius norm $\|Q^\tT P\|_F^2$ quantifies the total squared alignment between the two orthonormal frames, and thus serves as a measure of the closeness between the subspaces $\operatorname{span}(Q)$ and $\operatorname{span(P)}$. In particular, $\|Q^\tT P\|_F^2=k$ if and only if $\operatorname{span}(Q)=\operatorname{span}(P)$, and decreases as the principal angles between the subspaces increase.

In the numerical experiments reported below, we set $n=20$ and $k=5$, so that the ambient dimension is $d=nk=100$. Three of the test functions are nonconvex or multimodal, and are designed to probe the algorithm’s ability to escape local minima and efficiently explore the Grassmannian landscape. The implementation is written in Python, using NumPy version 2.2.4, and executed under the Windows Subsystem for Linux (WSL) environment.

\subsection{PCA objective function}

Let $\Sigma \in \mathrm{Mat}_n(\bR)$ be a symmetric positive semidefinite matrix, typically interpreted as a covariance matrix. We define the following function on the Grassmannian $\mathrm{Gr}(k,n)$:
\begin{equation}
f(V) = \tr(Q^\tT \Sigma Q),
\end{equation}
where $Q \in \mathrm{Mat}_{n,k}(\bR)$ is an orthonormal basis for the subspace $V$. The function $f(V)$ is minimized when $V$ spans the top $k$ eigenvectors of $\Sigma$, corresponding to the directions of greatest variance. At the global minimum, the objective value equals the negative of the sum of the top $k$ eigenvalues of $\Sigma$.

For the purpose of testing our optimization algorithm, we construct $\Sigma$ as a diagonal matrix with positive entries:
\begin{equation*}
\Sigma = \operatorname{diag}(20,19,\ldots,1).
\end{equation*}
The result produced by the algorithm is the matrix
\begin{equation}
Q^\ast =
\begin{bmatrix}
-0.0000 & 0.0000 & 0.0000 & 0.0000 & 0.0000 \\
-0.0000 & -0.0000 & 0.0000 & -0.0000 & 0.0000 \\
-0.0000 & -0.0001 & -0.0000 & -0.0000 & 0.0000 \\
-0.0000 & 0.0000 & 0.0000 & -0.0000 & 0.0001 \\
-0.0000 & -0.0000 & -0.0000 & 0.0000 & -0.0000 \\
0.0000 & 0.0000 & -0.0000 & 0.0000 & 0.0000 \\
0.0000 & -0.0000 & 0.0001 & -0.0000 & -0.0000 \\
0.0000 & -0.0001 & -0.0001 & 0.0000 & -0.0000 \\
-0.0000 & 0.0000 & -0.0001 & 0.0000 & 0.0000 \\
-0.0000 & 0.0001 & 0.0000 & 0.0001 & -0.0000 \\
0.0000 & 0.0000 & 0.0000 & -0.0000 & 0.0000 \\
0.0000 & 0.0000 & -0.0000 & 0.0000 & -0.0000 \\
-0.0000 & -0.0000 & 0.0000 & 0.0000 & -0.0001 \\
0.0000 & -0.0000 & -0.0000 & -0.0001 & 0.0001 \\
0.0000 & -0.0000 & 0.0001 & 0.0000 & 0.0000 \\
0.0157 & -0.0368 & 0.9893 & -0.1234 & 0.0665 \\
0.0516 & -0.0162 & -0.0692 & -0.0059 & 0.9961 \\
-0.1014 & 0.9913 & 0.0450 & 0.0660 & 0.0249 \\
0.0796 & -0.0632 & 0.1189 & 0.9876 & 0.0090 \\
-0.9902 & -0.1080 & 0.0170 & 0.0704 & 0.0511
\end{bmatrix}
\end{equation}
(with entries rounded to four decimal places), and the corresponding optimal value is
\begin{equation}
f(V^\ast)=89.99999931,
\end{equation}
which is in excellent agreement with the theoretical maximum of $90$.

\subsection{Cordal distance to a fixed subspace}

Let $P \in \mathrm{Mat}_{n,k}(\bR)$ be an orthonormal matrix whose columns span a fixed reference subspace $V_{ref} \in \mathrm{Gr}(k,n)$. Given a variable subspace $V = \operatorname{span}(Q)$, where $Q \in \mathrm{Mat}_{n,k}(\bR)$ is orthonormal, we define the objective function as the squared chordal distance \cite{TVSC11}:
\begin{equation}
\begin{split}
f(V) &= k - \| Q^\tT P \|_F^2\\
& =\tr(I_k - Q^\tT P P^\tT Q).
\end{split}
\end{equation}
This function is a smooth, rotation invariant measure of the angular discrepancy between the subspaces $V$ and $V_{ref}$. It achieves its minimum value of zero when $V = V_{ref}$, i.e., when $\operatorname{span}(Q)=\operatorname{span}(P)$. Geometrically, it corresponds to the squared chordal distance derived from the principal angles between the two subspaces.

For the purpose of testing our optimization algorithm, we fix a reference subspace $V_{ref} = \operatorname{span}(P)$ by randomly drawing $P=P_1\in \mathrm{Mat}_{n,k}(\bR)$ from the Stiefel manifold $\mathrm{V}_k(\bR^n)$:
\begin{equation}\label{p1_def}
P_1 =
\begin{bmatrix}
-0.1502 & 0.1246 & -0.1858 & 0.2688 & -0.2283 \\
0.2359 & -0.0589 & 0.2190 & 0.2316 & -0.2573 \\
-0.1873 & 0.2200 & -0.2186 & 0.0325 & 0.2314 \\
0.2177 & -0.2193 & 0.1236 & -0.1784 & -0.2835 \\
-0.2634 & -0.2732 & 0.1612 & -0.2198 & -0.3349 \\
-0.2404 & -0.3304 & -0.1157 & -0.2327 & 0.0783 \\
0.3175 & -0.0554 & -0.2848 & -0.1629 & 0.1717 \\
-0.2523 & -0.3586 & -0.0888 & -0.0295 & 0.3059 \\
-0.2005 & 0.2531 & -0.0529 & 0.1442 & -0.1602 \\
0.0254 & -0.1304 & 0.1385 & -0.2107 & 0.1724 \\
0.1812 & -0.1600 & 0.2000 & -0.1225 & -0.1943 \\
0.0416 & -0.3669 & 0.1466 & 0.1045 & -0.2294 \\
-0.3154 & 0.1979 & -0.1573 & -0.3159 & -0.3247 \\
0.0743 & -0.2198 & -0.1068 & -0.1548 & 0.2115 \\
0.3148 & -0.1081 & -0.2910 & 0.1124 & -0.0073 \\
0.0779 & -0.1567 & -0.2197 & 0.5016 & -0.1341 \\
0.2904 & -0.1657 & -0.5193 & -0.0795 & -0.1225 \\
-0.2637 & -0.1801 & -0.3949 & -0.0870 & -0.1839 \\
-0.1490 & -0.2496 & 0.2044 & 0.3056 & 0.3657 \\
0.2933 & 0.2863 & 0.0780 & -0.3535 & 0.0994
\end{bmatrix}.
\end{equation}
We then attempt to recover this subspace by minimizing $f(V)$ over the Grassmannian $\mathrm{Gr}(k,n)$. The optimization algorithm returns the solution:
\begin{equation}
Q^\ast =
\begin{bmatrix}
-0.1197 & 0.0630 & -0.4011 & -0.0031 & 0.1314 \\
-0.4285 & -0.1569 & 0.0461 & 0.0275 & 0.1277 \\
0.3349 & 0.1653 & -0.1500 & 0.0013 & 0.1543 \\
-0.3316 & -0.0687 & 0.1668 & -0.1080 & -0.2634 \\
-0.1785 & -0.0550 & -0.0988 & 0.1953 & -0.4951 \\
0.2463 & -0.1402 & -0.0613 & -0.0285 & -0.3942 \\
0.1408 & 0.0287 & 0.1823 & -0.4330 & 0.0119 \\
0.3902 & -0.3151 & -0.0595 & 0.0236 & -0.1936 \\
-0.0596 & 0.1956 & -0.2788 & 0.1423 & 0.1196 \\
0.1238 & -0.0695 & 0.2546 & 0.0476 & -0.1525 \\
-0.2779 & -0.0758 & 0.1933 & -0.0077 & -0.1748 \\
-0.2785 & -0.3317 & -0.0394 & 0.0172 & -0.1841 \\
-0.0274 & 0.4391 & -0.2527 & 0.0685 & -0.3298 \\
0.2085 & -0.1364 & 0.1431 & -0.1797 & -0.1332 \\
-0.0594 & -0.1074 & -0.0183 & -0.4233 & 0.1194 \\
-0.1644 & -0.3088 & -0.3554 & -0.1941 & 0.2491 \\
-0.0326 & 0.0072 & -0.1115 & -0.6163 & -0.0982 \\
0.1036 & 0.0281 & -0.3969 & -0.1908 & -0.3102 \\
0.2346 & -0.4640 & 0.0236 & 0.2412 & 0.1520 \\
0.0062 & 0.3512 & 0.4164 & -0.0973 & 0.0352
\end{bmatrix}
\end{equation}
(with entries rounded to four decimal places). The corresponding optimal function value is
\begin{equation}
f(V^\ast)=0.00016706,
\end{equation}
which is in excellent agreement with the theoretical minimum of $0$. Although $Q^\ast$ is not numerically close to $P_1$, it spans the same subspace in $ \mathrm{Gr}(k,n)$. To verify this, we compute the alignment matrix $R=(P^\tT_1 Q)^{-1}$, and confirm that it is approximately orthogonal with $\|R^\tT R-I_k\|_F=8.28\times 10^{-5}$.

\subsection{Bimodal subspace alignment function}

Let $P_1, P_2 \in \mathrm{Gr}(k,n)$ be two fixed and distinct subspaces, and let $V$ be a variable subspace. We define the bimodal alignment objective as
\begin{equation}
f(V) = \max\big\{ \| Q^\tT P_1 \|_F^2,\ \| Q^\tT P_2 \|_F^2 \big\}.
\end{equation}
This function measures the strongest total squared alignment of $V$ with either of the two target subspaces. The function $f(V)$ typically exhibits local minima in neighborhoods of both $P_1$ and $P_2$. 

To test our optimization algorithm, we randomly generate two subspaces $P_1, P_2 \in \mathrm{Gr}(k,n)$. Specifically, let $P_1$ be the (previously selected) matrix given by \eqref{p1_def}, and define
\begin{equation}\label{p2_def}
P_2 =
\begin{bmatrix}
-0.3143 & 0.2694 & -0.3297 & -0.0142 & -0.0316 \\
0.2495 & -0.4235 & -0.1919 & 0.0592 & 0.2191 \\
0.1312 & -0.2700 & -0.1111 & -0.2714 & 0.3331 \\
0.0052 & 0.0628 & 0.3018 & -0.0509 & 0.2455 \\
-0.0178 & 0.2174 & -0.4209 & -0.0639 & 0.0344 \\
-0.1033 & -0.2140 & -0.0939 & 0.0698 & 0.1463 \\
-0.3831 & -0.3233 & 0.0490 & 0.2698 & -0.1727 \\
0.1377 & -0.1426 & -0.2619 & -0.2988 & -0.1954 \\
-0.0264 & -0.0658 & -0.2540 & 0.2472 & 0.1171 \\
0.1635 & -0.2928 & 0.1476 & -0.2246 & 0.2523 \\
-0.0590 & -0.4255 & -0.2355 & 0.3699 & -0.1399 \\
-0.1270 & -0.2882 & 0.0749 & -0.3339 & -0.3113 \\
-0.1377 & -0.1466 & 0.4141 & -0.1233 & -0.1324 \\
0.1962 & -0.0535 & 0.0452 & 0.3463 & -0.4008 \\
0.2896 & 0.2502 & 0.0689 & 0.0445 & -0.1669 \\
-0.2139 & -0.0713 & 0.0357 & -0.2400 & -0.1482 \\
-0.1248 & 0.0012 & -0.0471 & 0.1948 & 0.2094 \\
0.4086 & 0.0366 & -0.1849 & 0.0221 & 0.0743 \\
-0.3234 & 0.0763 & 0.1542 & 0.1983 & 0.4666 \\
-0.3518 & 0.0003 & -0.3270 & -0.3416 & -0.0133
\end{bmatrix}.
\end{equation}
Evaluating the objective function at the two reference subspaces gives
\begin{equation}
\begin{split}
f(P_1) &= f(P_2)\\
&= 8.41823649.
\end{split}
\end{equation}

The algorithm converges to the solution
\begin{equation}
Q^\ast =
\begin{bmatrix}
-0.0145 & -0.4123 & -0.0744 & -0.2759 & -0.0527 \\
0.3464 & -0.0824 & -0.5397 & 0.0935 & -0.1163 \\
0.3073 & 0.3442 & -0.2346 & -0.4228 & -0.0959 \\
0.1282 & 0.5239 & 0.1473 & -0.1272 & 0.3165 \\
0.0698 & -0.1352 & 0.4090 & -0.3758 & -0.1368 \\
-0.1015 & -0.0726 & -0.1962 & -0.1782 & -0.1319 \\
0.0443 & -0.0444 & -0.0866 & -0.1270 & 0.0698 \\
0.4411 & -0.2630 & 0.0258 & 0.1082 & 0.1536 \\
0.1154 & -0.1640 & -0.0304 & -0.1722 & 0.1526 \\
0.2398 & -0.0087 & -0.1624 & -0.2780 & -0.0635 \\
0.0575 & -0.3008 & 0.0695 & 0.1532 & -0.1730 \\
-0.0532 & 0.0105 & -0.3952 & -0.1715 & 0.1251 \\
0.0101 & 0.0815 & -0.1108 & 0.1366 & -0.2075 \\
-0.0525 & 0.1117 & -0.0665 & 0.1520 & -0.6003 \\
0.1275 & -0.2402 & 0.1973 & 0.0463 & -0.1169 \\
0.3322 & 0.2625 & 0.3523 & -0.0584 & -0.2633 \\
-0.3330 & -0.0892 & -0.0895 & -0.2589 & 0.1822 \\
0.2253 & 0.0618 & -0.0835 & 0.4803 & 0.3377 \\
-0.1635 & -0.0584 & 0.0885 & -0.0660 & 0.2707 \\
0.4029 & -0.2388 & 0.1475 & -0.1010 & 0.1732
\end{bmatrix}.
\end{equation}
For this solution, we find:
\begin{equation}
\begin{split}
\| Q^{\ast\tT} P_1 \|_F^2 &= 0.02126794, \\
\| Q^{\ast\tT} P_2 \|_F^2 &= 1.48077098,
\end{split}
\end{equation}
and hence $f(Q^\ast) = 1.48077098$. Although the function $f(V)$ is symmetric in $P_1$ and $P_2$, the global minimum is closer to $P_1$ than $P_2$.

\subsection{Log-determinant volume penalty}

Let $A \in \mathrm{Mat}_n(\bR)$ be a symmetric positive definite matrix. We define the following objective function on the Grassmannian:
\begin{equation}
f(V) = -\log \det(Q^\tT A Q),
\end{equation}
where $Q \in \mathrm{Mat}_{n,k}(\bR)$ is an orthonormal basis for the subspace $V \in \mathrm{Gr}(k,n)$.

This function is closely related to the volume of the ellipsoid obtained by projecting the positive definite matrix $A$ onto the subspace $V$. Specifically, $\det(Q^\tT A Q)$ is proportional to the squared volume of the image of the unit ball in $V$ under the linear map induced by $A^{1/2}$. Therefore, minimizing $f(V)$ amounts to identifying the $k$-dimensional subspace that preserves the largest volume under projection by $A$.

While the function is convex in the projection matrix $QQ^\tT$, it is generally nonconvex over the Grassmannian $\mathrm{Gr}(k,n)$ due to its quotient geometry. When the eigenvalues of $A$ are tightly clustered or highly disparate, the landscape of $f$ may feature multiple local minima, especially when several subspaces correspond to nearly equivalent projected volumes.

To evaluate the performance of our optimization algorithm on this objective, we construct $A$ with a synthetic spectrum designed to induce multiple near-optimal subspaces:
\begin{equation}
A = \operatorname{diag}(10, 9, 1, \dots, 1).
\end{equation}
This configuration creates a dominant eigenspace spanned by the top eigenvectors and tests the algorithm’s ability to locate the most volume-preserving subspace in the presence of curvature and eigenvalue degeneracy.

The optimizer returns the following solution:
\begin{equation}
Q^\ast =
\begin{bmatrix}
-0.0000 & 1.0000 & 0.0000 & 0.0000 & -0.0000 \\
-1.0000 & -0.0000 & 0.0000 & 0.0000 & -0.0000 \\
-0.0000 & 0.0000 & -0.0006 & -0.6978 & 0.1248 \\
0.0000 & 0.0000 & -0.2449 & -0.2658 & -0.5333 \\
-0.0000 & -0.0000 & -0.2383 & 0.0778 & -0.1336 \\
-0.0000 & -0.0000 & 0.4749 & -0.1083 & -0.0164 \\
-0.0000 & -0.0000 & 0.1160 & -0.0313 & -0.2383 \\
0.0000 & -0.0000 & 0.1045 & 0.1411 & 0.0541 \\
-0.0000 & 0.0000 & -0.4622 & 0.0327 & 0.4064 \\
0.0000 & 0.0000 & -0.0516 & 0.1181 & -0.2424 \\
0.0000 & -0.0000 & -0.1200 & -0.1575 & -0.2613 \\
-0.0000 & 0.0000 & 0.1512 & 0.2889 & 0.0140 \\
0.0000 & -0.0000 & -0.2514 & 0.2499 & -0.3951 \\
-0.0000 & -0.0000 & -0.2896 & -0.0372 & -0.1839 \\
-0.0000 & 0.0000 & -0.0973 & 0.2897 & 0.0617 \\
0.0000 & 0.0000 & -0.1708 & 0.1217 & -0.0481 \\
0.0000 & 0.0000 & -0.1419 & -0.2601 & -0.0222 \\
0.0000 & 0.0000 & 0.2761 & 0.1482 & -0.2051 \\
-0.0000 & 0.0000 & 0.2751 & 0.0003 & -0.2910 \\
0.0000 & 0.0000 & -0.1481 & 0.1680 & 0.0808
\end{bmatrix}
\end{equation}
The corresponding optimal value is:
\begin{equation}
f(Q^\ast) = -4.4998,
\end{equation}
which is consistent with the theoretical lower bound determined by the log-volume of the dominant eigen-directions.

\subsection{Subspace clustering residual}

Let $X_1, \dots, X_r \in \mathrm{Mat}_{n,m}(\bR)$ be a collection of data matrices, each consisting of $m$ data vectors drawn from a distinct low-dimensional subspace of $\bR^n$. Let $V = \operatorname{span}(Q) \in \mathrm{Gr}(k,n)$ be a candidate subspace, represented by an orthonormal matrix $Q \in \mathrm{Mat}_{n,k}(\bR)$. Given a fixed data matrix $X \in \mathrm{Mat}_{n,m}(\bR)$, we define the subspace clustering residual function on the Grassmannian as:
\begin{equation}
f_X(V) = \min_{1 \leq i \leq r} \| X - QQ^\tT X_i \|_F^2.
\end{equation}
Here, $QQ^\tT$ is the orthogonal projection onto the subspace $V$, and $f_X(V)$ measures the residual error of projecting the candidate data matrix $X$ onto each of the reference subspaces $\operatorname{span}(X_i)$, returning the smallest among these projection errors.

This objective function is piecewise smooth: it is differentiable within regions where the minimizing index $i$ remains constant, but nonsmooth at the boundaries separating these regions. Each subspace $\operatorname{span}(X_i)$ typically induces a local minimum of the function, leading to a multimodal optimization landscape consisting of sharp attraction basins - each centered around one of the $X_i$ - separated by nonsmooth ridges. This structure closely models the geometry of subspace clustering problems, where the goal is to identify the best-matching subspace among a finite set of candidates.

To evaluate the performance of our optimization algorithm in this setting, we consider the case $r = 3$ and $m = 12$, and construct synthetic data matrices $X_i = P_i Z_i$, where $P_i \in \mathrm{V}_k(\bR^n)$ are fixed orthonormal basis matrices representing distinct subspaces, and $Z_i \in \mathrm{Mat}_{k,m}(\bR)$ are matrices with i.i.d. standard Gaussian entries. The matrices $P_1$ and $P_2$ are as defined in \eqref{p1_def} and \eqref{p2_def}, while $P_3$ is given by:
\begin{equation}
P_3 =
\begin{bmatrix}
-0.1390 & -0.1537 & -0.3595 &  0.1077 & -0.1327 \\
0.3579 & -0.3058 & -0.1719 &  0.2780 & -0.1325 \\
0.1792 & -0.0517 & -0.0835 & -0.0192 & -0.2234 \\
-0.1217 &  0.0481 & -0.0312 & -0.0190 & -0.2553 \\
0.0480 &  0.0721 &  0.3349 & -0.1425 & -0.5004 \\
0.0613 &  0.1127 & -0.2846 &  0.1409 & -0.1942 \\
-0.3438 & -0.2623 &  0.0629 &  0.3949 & -0.2553 \\
0.0982 & -0.3549 &  0.3214 & -0.1534 &  0.0771 \\
0.2514 &  0.3808 &  0.4924 &  0.1946 & -0.1096 \\
-0.0955 &  0.0599 & -0.0881 & -0.3908 & -0.1217 \\
0.1858 & -0.0130 & -0.2281 & -0.2885 &  0.3168 \\
0.3453 &  0.1694 & -0.1641 &  0.3670 &  0.2529 \\
0.2726 & -0.2969 & -0.0442 & -0.1964 & -0.2372 \\
-0.0720 &  0.2829 & -0.0233 &  0.0393 & -0.2494 \\
0.3083 & -0.2587 & -0.0587 &  0.0192 & -0.1084 \\
0.0901 &  0.1662 & -0.1897 & -0.3868 & -0.1561 \\
0.2764 & -0.1413 &  0.3400 &  0.0447 &  0.1874 \\
0.1009 &  0.2984 & -0.1802 &  0.2640 & -0.1090 \\
0.3156 & -0.0803 & -0.0893 & -0.0274 & -0.3126 \\
0.2773 &  0.3309 & -0.0759 & -0.1401 &  0.0012
\end{bmatrix}.
\end{equation}
Each objective function $f_{X_i}(V)$ achieves its minimum at the corresponding ground truth subspace $V =\operatorname{span}(P_i)$, i.e.,
\begin{equation}\label{minfxi}
f_{X_i}(X_i) = 0.
\end{equation}

This experimental setup provides a rigorous test of the algorithm’s ability to recover subspace structure in the presence of multiple competing models and nonsmooth landscape features.

The optimization algorithm successfully recovers each ground truth subspace, yielding the following results:
\begin{equation*}
\begin{split}
\mathop{\arg\min} f_{X_1}(V) = V_1^\ast, &\quad \text{with } f_{X_1}(V_1^\ast) = 1.06 \times 10^{-6}, \\
\mathop{\arg\min} f_{X_2}(V) = V_2^\ast, &\quad \text{with } f_{X_2}(V_2^\ast) = 7.00 \times 10^{-7}, \\
\mathop{\arg\min} f_{X_3}(V) = V_3^\ast, &\quad \text{with } f_{X_3}(V_3^\ast) = 5.09 \times 10^{-7},
\end{split}
\end{equation*}
demonstrating excellent agreement with the theoretical minima \eqref{minfxi}, and confirming the algorithm’s ability to identify the correct subspace in each case.


\begin{thebibliography}{99}
\bibitem{AMS08} Absil, P. -A., Mahony, R., and Sepulchre, R.: \textit{Optimization Algorithms on Matrix Manifolds}, Princeton University Press (2008).
\bibitem{BZA24} Bendokat, T., Zimmermann, R., and Absil, P.-A.: A Grassmann manifold handbook: basic geometry and computational aspects, \textit{Adv. Comp. Math.}, \textbf{50:6} (2024).
\bibitem{B23} Boumal, N.: \textit{An Introduction to Optimization on Smooth Manifolds}, Cambridge University Press (2023).
\bibitem{BGBMZ06} Brest, J., Greiner, S., Boskovic, B., Mernik, M., and Zumer, V.: Self-Adapting Control Parameters in Differential Evolution: A Comparative Study on Numerical Benchmark Problems, \textit{IEEE Trans. Evolutionary Comp.}, \textbf{10}, 646 - 657 (2006).
\bibitem{EV13} Elhamifar, E., and Vidal, R.: Sparse Subspace Clustering: Algorithm, Theory, and Applications, \textit{IEEE Trans, on Pattern Analysis and Machine Intelligence}, \textbf{35}, 2765 - 2781 (2013).
\bibitem{FWLLC23} Fei, Y., Wei, X., Liu, Y., Li, Z., and Chen, M.: A Survey of Geometric Optimization for Deep Learning: From Euclidean Space to Riemannian Manifold, arXiv:2302.08210v1 (2023).
\bibitem{SP97} Storn, R., and Price, K.: Differential evolution - a simple and efficient heuristic for global optimization over continuous spaces, \textit{J. Global Optim.}, \textbf{11}, 341 - 359 (1997).
\bibitem{TVSC11} Turaga, P., Veeraraghavan, A., Srivastava, A., and Chellappa, R.: Statistical Computations on Grassmann and Stiefel Manifolds for Image and Video-Based Recognition, \textit{IEEE Trans. on Pattern Analysis and Machine Intelligence}, \textbf{33}, 2273 - 2286, (2011).
\bibitem{ZZHH18} Zhang, J., Zhu, G., Heath Jr., R. W., and Huang, K.: Grassmannian Learning: Embedding Geometry Awareness in Shallow and Deep Learning, arXiv:1808.02229v2 (2018).
\end{thebibliography}
\end{document}